\documentclass[12pt]{amsart}

\usepackage{amsmath}
\usepackage{amsthm}
\usepackage{amssymb}
\usepackage{graphics} 
\usepackage{pgf}
\usepackage{tikz}
\usepackage{pifont}
\usetikzlibrary{arrows}
\usepackage[all]{xy}
\usepackage{lineno}

\def\frk{\frak}               

\def\Phi{{\frk n}}
\def\Phi{{\frk N}}
\def\opn#1#2{\def#1{\operatorname{#2}}} 
\opn\chara{char} \opn\length{\ell} \opn\pd{pd} \opn\rk{rk}
\opn\projdim{proj\,dim} \opn\injdim{inj\,dim} \opn\rank{rank}
\opn\depth{depth} \opn\grade{grade} \opn\height{height}
\opn\embdim{emb\,dim} \opn\codim{codim}

\opn\Tr{Tr} \opn\bigrank{big\,rank}
\opn\superheight{superheight}\opn\lcm{lcm}
\opn\trdeg{tr\,deg}
\opn\reg{reg} \opn\lreg{lreg} \opn\ini{in} \opn\lpd{lpd}
\opn\size{size}\opn\bigsize{bigsize}
\opn\cosize{cosize}\opn\bigcosize{bigcosize}
\opn\sdepth{sdepth}\opn\sreg{sreg}
\opn\link{link}\opn\fdepth{fdepth}
\opn\index{index}
\opn\index{index}
\opn\indeg{indeg}
\opn\N{N}
\opn\mult{mult}
\opn\SSC{SSC}
\opn\SC{SC}
\opn\lk{lk}
\opn\HS{HS}
\opn\div{div} \opn\Div{Div} \opn\cl{cl} \opn\Cl{Cl}
\opn\Spec{Spec} \opn\Supp{Supp} \opn\supp{supp} \opn\Sing{Sing}
\opn\Ass{Ass} \opn\Min{Min}\opn\Mon{Mon} \opn\dstab{dstab} \opn\astab{astab}
\opn\Syz{Syz}
\opn\reg{reg}
\opn\Ann{Ann} \opn\Rad{Rad} \opn\Soc{Soc}
\opn\Im{Im} \opn\Ker{Ker} \opn\Coker{Coker} \opn\Am{Am}
\opn\Hom{Hom} \opn\Tor{Tor} \opn\Ext{Ext} \opn\End{End}\opn\Der{Der}
\opn\Aut{Aut} \opn\id{id}

\opn\nat{nat}
\opn\pff{pf}
\opn\Pf{Pf} \opn\GL{GL} \opn\SL{SL} \opn\mod{mod} \opn\ord{ord}
\opn\Gin{Gin} \opn\Hilb{Hilb}\opn\sort{sort}
\opn\initial{init}
\opn\ende{end}
\opn\height{height}
\opn\type{type}
\opn\aff{aff} \opn\con{conv} \opn\relint{relint} \opn\st{st}
\opn\lk{lk} \opn\cn{cn} \opn\core{core} \opn\vol{vol}
\opn\link{link} \opn\Link{Link}\opn\lex{lex}
\opn\gr{gr}

\def\pot#1#2{#1[\kern-0.28ex[#2]\kern-0.28ex]}

\opn\dirlim{\underrightarrow{\lim}}
\opn\inivlim{\underleftarrow{\lim}}
\let\to=\rightarrow

\def\Implies{\ifmmode\Longrightarrow \else
        \unskip${}\Longrightarrow{}$\ignorespaces\fi}
\def\implies{\ifmmode\Rightarrow \else
        \unskip${}\Rightarrow{}$\ignorespaces\fi}
\def\iff{\ifmmode\Longleftrightarrow \else
        \unskip${}\Longleftrightarrow{}$\ignorespaces\fi}

\let\:=\colon
\newtheorem{Theorem}{Theorem}[section]
 \newtheorem{Lemma}[Theorem]{Lemma}
 \newtheorem{Corollary}[Theorem]{Corollary}
 \newtheorem{Proposition}[Theorem]{Proposition}
 \newtheorem{Remark}[Theorem]{Remark}

 \newtheorem{Definition}[Theorem]{Definition}
  \newtheorem{Question}[Theorem]{Question}

\newtheorem{Notation}[Theorem]{Notation}
 \newtheorem{Conventions}[Theorem]{Conventions}
\let\epsilon\varepsilon
\let\kappa=\varkappa
\def\qed{\ifhmode\textqed\fi
      \ifmmode\ifinner\quad\qedsymbol\else\dispqed\fi\fi}
\def\textqed{\unskip\nobreak\penalty50
       \hskip2em\hbox{}\nobreak\hfil\qedsymbol
       \parfillskip=0pt \finalhyphendemerits=0}
\def\dispqed{\rlap{\qquad\qedsymbol}}

\opn\dis{dis}
\def\pnt{{\raise0.5mm\hbox{\large\bf.}}}

\opn\Lex{Lex}

\opn\Spec{Spec} \opn\Supp{Supp} \opn\supp{supp}
 \opn\Ass{Ass}
 \opn\p{Ass}
 \opn\min{min}
 \opn\max{max}
 \opn\MIN{Min}
 \opn\p{\mathfrak{p}}
\opn\Deg{Deg}

\begin{document}

 \title{Betti Numbers of Edge Ideals of Weighted Oriented Crown Graphs}

\author {Zexin Wang}

\address{School  of Mathematical Science, Soochow University, 215006 Suzhou, P.R.China}
\email{zexinwang6@outlook.com}

 \author{Dancheng Lu}

\address{School  of Mathematical Sciences, Soochow University, 215006 Suzhou, P.R.China}
\email{ludancheng@suda.edu.cn}

 \begin{abstract}We compute the multigraded Betti numbers of edge ideals for weighted oriented crown graphs using a novel approach. This approach, which we still call the \emph{induced subgraph approach}, originates from our prior work on computing Betti numbers of normal edge rings(see \cite{WL}). Notably, we prove that the total Betti numbers of edge ideals for weighted oriented crown graphs are independent of the weight function.
 \end{abstract}



\maketitle

\section{Introduction}
A weighted oriented graph is a triple $D = (V(D), E(D), w)$ where $V(D) = \{ x_1, \ldots, x_n \}$ is the set of vertices, $E(D)$ is a directed edge set, and $w$ is a weight function $w: V(D) \rightarrow \mathbb{N}^+$. Specifically, $E(D)$ consists of ordered pairs $(x_i, x_j) \in V(D) \times V(D)$  where the pair $(x_i, x_j)$ represents a directed edge from the vertex $x_i$ to the vertex $x_j$. We refer to the graph $G$ obtained by removing the orientation and weights from a weighted oriented graph $D$ as the \emph{underlying graph} of $D$. In this paper, we consider only weighted oriented simple graphs, meaning that their underlying graph is a simple graph (recall that a simple graph is a finite graph without loops or multiple edges). Let $R:=\mathbb{K}[x_1,\ldots,x_n]$ be a  standard graded polynomial ring over a field $\mathbb{K}$.
The edge ideal of a weighted oriented graph $D$ is a monomial ideal given by
\[
I(D) = (x_i x_j^{w(x_j)} \mid (x_i, x_j) \in E(D)) \subseteq R.
\]
If $w(x_j) = 1$ for all $j$, then $I(D)$ is the edge ideal of an unweighted undirected graph, which was introduced by Villarreal \cite{V} and has been extensively studied.

Edge ideals of weighted oriented graphs arose in the theory of Reed-Muller codes as initial ideals of vanishing ideals of projective spaces over finite fields(\cite{CNL},\cite{MPV2}). The edge ideal of a weighted oriented graph is a relatively new object of study in combinatorial commutative algebra, first investigated by Gimenez-Bernal-Simis-Villarreal-Vivares in \cite{GMSVV} in 2018. Subsequently, several papers further explored its algebraic properties and invariants, such as \cite{CK}, \cite{HLMRV}, \cite{MP}, \cite{PRT}, \cite{WZXZ}, and \cite{ZXWT}.

This paper focuses on the graded and multigraded Betti numbers of edge ideals of weighted oriented graphs. The seminal work of Hochster \cite{H} in 1977 pioneered the use of simplicial homology to compute multigraded Betti numbers of square-free monomial ideals, which led to the celebrated Hochster's formula for computing multigraded Betti numbers of arbitrary monomial ideals. Since edge ideals of simple graphs are precisely quadratic square-free monomial ideals, the Hochster's formula admits significant simplification in this case (see, e.g., \cite[Section 2.2]{EO}). For this reason, many studies have utilized Hochster's formula to investigate the Betti numbers of edge ideals of certain (unweighted, undirected) graphs and obtained significant results, such as \cite{EO}, \cite{R}, \cite{RS1}, and \cite{RS2}.

However, edge ideals of general weighted oriented graphs lose the square-free property, making Hochster's formula less effective for deriving explicit formulas for their Betti numbers. The Betti splitting technique, first introduced by Eliahou and Kervaire \cite{EK} for computing Betti numbers of monomial ideals, was later adapted by H\`{a} and Tuyl \cite{HT2} to edge ideals of (unweighted, unoriented) graphs. Their approach transforms the computation of Betti numbers into combinatorial problems on associated simple graphs, thereby circumventing the highly nontrivial task of computing homology group dimensions. Casiaday and Kara \cite{CK} extended this methodology to weighted oriented graphs, obtaining formulas for projective dimension, regularity, and recursive relations for graded Betti numbers of certain classes of these graphs.

While Betti splitting remains an effective computational tool that we will extensively employ, the complexity introduced by weights in weighted oriented graphs poses significant challenges. Even when recursive formulas for graded Betti numbers are obtainable, deriving complete and explicit general formulas remains difficult. Moreover, as discussed in \cite[Section 4]{FHT}, while Betti splittings may be prevalent among monomial ideals, proving their existence theoretically is often challenging (though sufficient conditions have been established in works such as \cite{B} and \cite{FHT}). To our knowledge, these obstacles explain why no non-trivial class of weighted oriented graphs has had its complete graded and multigraded Betti numbers fully characterized in the literature.

To address these challenges, we develop a novel approach that combines existing methods with our induced-subgraph technique originally devised for computing normal edge rings (see \cite{WL}, \cite{WL1}). We maintain the nomenclature \emph{induced-subgraph approach} for this unified framework, which offers new insights into the multigraded Betti numbers of edge ideals of weighted oriented graphs.

\begin{Lemma}\em Let $H$ be  an induced subgraph of a weighted oriented graph $D$.
Then, for $\mathbf{a} \in \mathbb{N}^{V(D)}$ with $\mathrm{supp}(\mathbf{a}) \subseteq V(H)$, where $\mathrm{supp}(\mathbf{a})$ is defined as $\mathrm{supp}(x^{\mathbf{a}})$, we have
\[
\beta_{i,\mathbf{a}}(I(H)) = \beta_{i,\mathbf{a}}(I(D)) \quad \text{for all } i \geq 0.
\]
\end{Lemma}

To facilitate the presentation of our approach, we first introduce some additional definitions and notations. Let $D$ be a weighted oriented graph with $V(D) = \{x_1, \ldots, x_n\}$, and let $\mathrm{pdim}(I(D)) = p$. Subsequently, $\beta_p^{R}(I)$ is referred to as the {\it top total Betti number} of $I$, while $\beta_{p,\mathbf{a}}^{R}(I)$ for $\mathbf{a}\in \mathbb{N}^n$ are designated as the {\it top multigraded Betti numbers} of $I$. An oriented graph $H$ is called \emph{induced subgraph} of $D$ if $V(H)\subseteq V(D)$ and for any $x_i, x_j\in V(H)$, $(x_i, x_j)$ is an edge in $H$ if and only if $(x_i, x_j)$ is an edge in $D$. Let $V(H) = W$, then we say that $H$ is the induced subgraph of $D$ on the set $W$, and in this case, denote $H$ as $D_W$.

 For any induced subgraph $H$ of $D$ (including $D$ itself), define the monomial set
\[
\mathcal{N}_H^k(D) := \{x^h \in R \mid \beta_{p,h}(I(D)) = k \text{ for some } h \in \mathbb{N}^n\}.
\]
Then it immediately follows that the top total Betti number of $I(D)$ satisfies
$\beta_p(I(D)) = \sum_k k \cdot |\mathcal{N}_H^k(D)|$.

The detailed steps of our approach is as follows.
  \begin{enumerate}
 \item Initially, first compute the projective dimension $p$ of $I(D)$, then obtain a suitable upper bound for the total Betti number $\beta_i(I(D))$.

  \item Secondly, compute the projective dimension of the edge ideal $I(H)$ and the sets $\mathcal{N}_H^k(D)$ for some induced subgraphs $H$ of $D$ (including $D$ itself). Denote this collection of induced subgraphs by $\mathcal{H}$.
  \item Next, for indices $i$ with $0 < i \leq p$, the following inequality holds:
\[
\sum_k k \cdot \left|\bigcup_{\mathrm{pdim}(I(H))=i}^{H\in \mathcal{H}} \mathcal{N}_H^k(D)\right| \leq \beta_i(I(D)).
\]
If computation of the left-hand side equals the upper bound obtained in the first step, then we achieve the goal of determining all multigraded Betti numbers of $I(D)$. In this case, $\beta_{i,h}(I(D))$ equals the top multigraded Betti number $\beta_{i,h}(I(H))$, where $H$ is the induced subgraph of $D$ on $\mathrm{supp}(h)$. Moreover, if $\beta_{i,h}(I(D))\neq 0$, then the projective dimension of $I(H)$ is equal to $i$.
  \item Otherwise, we must consider either additional induced subgraphs of $D$, or the multigraded Betti numbers $\beta_{i,h}(I(H))$ with lower homological degrees for edge ideals of induced subgraphs $H$ of $D$.
\end{enumerate}
In this paper, we successfully characterize the multigraded Betti numbers of edge ideals for weighted oriented crown graphs $G_n= (V(G_n), E(G_n), w)$ using this approach. Here, $V(G_n) = \{x_1,\ldots,x_n,y_1,\ldots,y_n\}$, $E(G_n) = \Big\{(x_i, y_j) \mid 1 \leq i, j \leq n, i \neq j \Big\}$, and $w$ is a weight function $w: V(G_n) \rightarrow \mathbb{N}^+$.  The underlying graph of $G_n$ is precisely the $n$\emph{-crown graph} $\mathcal{C}_{n,n}$, which is why we call $G_n$ a weighted oriented crown graph. Recall that the $n$-crown graph is defined as the bipartite graph obtained by removing a perfect matching $\{x_i,y_i\}_{1 \leq i \leq n}$ from the complete bipartite graph on vertex set $\{x_1,\ldots,x_n\} \cup \{y_1,\ldots,y_n\}$. The algebraic properties of graph ideals of crown graphs have been extensively studied, see for example \cite{KKSS}, \cite{RS1}, and \cite{RS2}. Throughout this paper, we consistently denote $w(y_i)$ by $w_i$ for $1 \leq i \leq n$, and define $I_n = (x_i y_j^{w_j} \mid 1 \leq i, j \leq n, i \neq j)$ as the edge ideal of $G_n$.

The main results of this paper is as follows.

\begin{Theorem}[\emph{Proposition \ref{Gs betti} and Theorem~\ref{main2}}]
Let $G_n$ be the weighted oriented crown graph as described previously, and let $I_n$ be its edge ideal. For all $0 \leq i$, the $i$-th total Betti number is:
\[
\beta_i(I_n) = \sum_{k=2}^n (k-1)2^{i+3-2k}\binom{n}{k}\binom{n-k}{i+3-2k} + (2^{i+2}-2)\binom{n}{i+2}.
\]
In particular, the projective dimension of $I_n$ is $2n-3$. Moreover, the total Betti numbers of the weighted crown graph $G_n$ are independent of the vertex weights.
\end{Theorem}

\begin{Theorem}[\emph{Theorem~\ref{main} and Corollary \ref{crown reg}}]
The multigraded Betti numbers of $I_n$ are given by
\[
\beta_{i, \alpha}(I_n) = \begin{cases}
k-1, & \text{if } \alpha \in \mathcal{N}_i^k(G_n) \text{ for some }2\leq k\leq n; \\
1, & \text{if } \alpha \in \mathcal{M}_i(G_n); \\
0, & \text{otherwise.}
\end{cases}
\]
Furthermore, the graded Betti numbers of $I_n$ are
\[
\beta_{i,j}(I_n) =\sum_{k=2}^n \sum_{\substack{\alpha \in \mathcal{N}_i^k(G_n) \\ \deg(\alpha)=j}} (k-1) + |\{\alpha \in \mathcal{M}_i(G_n) : \deg(\alpha)=j\}|.
\]
In particular, the regularity of $I_n$ is $\sum_{i=1}^n w_i - n + 3$. Here, $\mathcal{N}_i^k(G_n)$ and $\mathcal{M}_i(G_n)$ are monomial sets completely determined by certain induced subgraphs of $G_n$, as detailed in Definition~\ref{topsupp}.
\end{Theorem}

In \cite{RS1} and \cite{RS2}, Rather and Singh primarily employed Hochster's formula to compute the graded Betti numbers of edge ideals of crown graphs. When all vertex weights equal 1 in weighted oriented crown graphs, their edge ideals precisely coincide with those of crown graphs. Therefore, as a corollary of Theorem~\ref{main}, we recover the main results of \cite{RS1} and \cite{RS2} while bypassing the computationally intensive determination of simplicial homology dimensions (see Corollary~\ref{cor}). However, since we also computed the multigraded Betti numbers, we can better explain the Betti number formula using the combinatorial data of crown graphs.

Furthermore, Alesandroni \cite{A2} investigated minimal free resolutions of arbitrary monomial ideals, demonstrating that the multigraded Betti numbers of any monomial ideal decompose as sums of those from some dominant ideals and some purely nondominant ideals. Theorem 4.4 in \cite{A1} establishes that the Taylor resolution of a dominant ideal already constitutes its minimal multigraded free resolution. Consequently, the study of minimal monomial resolutions reduces to analyzing purely nondominant ideals - a generally challenging problem. Following Definition 2.5 in \cite{A2}, one immediately observes that $I_n$ is purely nondominant. Therefore, investigating its multigraded Betti numbers provides valuable insights into minimal free resolutions of purely nondominant ideals.

The paper is structured as follows. In Section 2, we introduce basic concepts and technical lemmas frequently used throughout this paper. In Section 3, following the first step of our approach, we obtain an upper bound for the total Betti numbers of $I_n$ in recursive form. In Section 4, implementing the second and third steps of our approach, we compute the top multigraded Betti numbers of edge ideals for certain induced subgraphs of $G_n$ and combine these results with Section 3 to derive formulas for the total Betti numbers, multigraded Betti numbers, and graded Betti numbers of $I_n$. Moreover, we actually obtain all multigraded Betti numbers of edge ideals for every induced subgraph of $G_n$.

\section{Preliminaries}

Let $\mathbb{K}$ be a field, and let $R:=\mathbb{K}[x_1,\ldots,x_n]$ be the polynomial ring in variables $x_1,\ldots,x_n$.   For any monomial ideal $I$ of $R$, there exists the minimal multigraded free resolution of $I$ that has the form:
\begin{equation*}\label{free}0\rightarrow\underset{\mathbf{a}\in\mathbb{N}^n}{\bigoplus}R[-\mathbf{a}]^{\beta_{p,\mathbf{a}}^{R}(I)}\rightarrow \cdots \rightarrow \underset{\mathbf{a}\in\mathbb{N}^n}{\bigoplus}R[-\mathbf{a}]^{\beta_{1,\mathbf{a}}^{R}(I)}\rightarrow\underset{\mathbf{a}\in\mathbb{N}^n}{\bigoplus}R[-\mathbf{a}]^{\beta_{0,\mathbf{a}}^{R}(I)}\rightarrow I \rightarrow 0.\end{equation*}
Note that $R[-\mathbf{a}]$ is the cyclic free $R$-module generated in degree $\mathbf{a}$.
The {\it projective dimension} of $I$ are defined to be  $\mathrm{pdim}_{R}\,(I):={\mbox{max}}\,\{i\mid \beta_{i,\,\mathbf{a}}^{R}(I)\neq 0 \mbox{ for some } \mathbf{a} \}.$
The number $\beta_{i,\mathbf{a}}^{R}(I)={\rm{dim}}_{\mathbb{K}}\mathrm{Tor}_i^R(I,\mathbb{K})_\mathbf{a}$ is called the $(i,\mathbf{a})$-th {\it multigraded Betti number} of $I$ and $\beta_{i}^{R}(I):=\sum_{\mathbf{a}\in \mathbb{N}^n}\beta_{i,\mathbf{a}}^{R}(I)$ is called the $i$-th {\it total Betti number} of $I$.
Furthermore, for $i \in \mathbb{N}$ and $j \in \mathbb{N}$, the numbers
\[
\beta_{i,j}^R(I) := \sum_{\substack{\mathbf{a} \in \mathbb{N}^n \\ \sum_k a_k = j}} \beta_{i,\mathbf{a}}^R(I)
\]
are called the \emph{graded Betti numbers} of $I$. The \emph{Castelnuovo-Mumford regularity} of $R/I$ is defined as
\[
\reg_R (I) := \max \{j - i : \beta_{i,j}^R(I) \neq 0\}.
\]

The \emph{multigraded Betti numbers} of $R/I$ are defined as $\beta_{i,\mathbf{a}}^R(R/I) = \beta_{i-1,\mathbf{a}}^R(I)$ for $i \in \mathbb{N}$ and $\mathbf{a} \in \mathbb{N}^n$. Consequently, one has $\mathrm{pdim}_R (R/I) =\mathrm{pdim}_R (I) +1 $ and $\reg_R (R/I) =\reg_R (I)-1$.
Denote by $p$ the projective dimension of $I$. Subsequently, $\beta_p^{R}(I)$ (resp. $\beta_{p+1}^{R}(R/I)$) is referred to as the {\it top total Betti number} of $I$ (resp. $R/I$), while $\beta_{p,\mathbf{a}}^{R}(I)$ (resp. $\beta_{p+1,\mathbf{a}}^{R}(R/I)$)for $\mathbf{a}\in \mathbb{N}^n$ are designated as the {\it top multigraded Betti numbers} of $I$ (resp. $R/I$). Furthermore, for any monomial ideal $I$ in $R$, we let $G(I)$ denote its minimal set of monomial generators. If it is clear from the context in which ring we are working, we will omit $R$ from the notation of the Betti numbers, projective dimension and regularity, and just write $\beta_{i,\mathbf{a}}(I)$, $\beta_i(I)$, $\beta_{i,j}(I)$, $\mathrm{pdim} (I)$ and $\reg (I)$. When we talk about a monomial $x^{\mathbf{a}}$ of $R$, we always mean $\mathbf{a} = (a_1, \ldots, a_n) \in \mathbb{N}^n$ and $x^{\mathbf{a}}=x_1^{a_1}\ldots x_n^{a_n}$.

\begin{Definition}\cite[Definition 1.1]{FHT}\label{betti splitting de}
Let $I, J,$ and $K$ be monomial ideals of $R$ such that $G(I)$ is the disjoint union of $G(J)$ and $G(K)$.
Then $I = J+K$ is a {\bf Betti splitting} if
\[\beta_{i,\mathbf{a}}(I) = \beta_{i,\mathbf{a}}(J)+\beta_{i,\mathbf{a}}(K)+\beta_{i-1,\mathbf{a}}(J \cap K) \hspace{.5cm}~~\mbox{for all $i\in \N$
and}\; \mathbf{a}\in\mathbb{N}^n.\]
\end{Definition}
When $I$ is a Betti splitting ideal, Definition \ref{betti splitting de} implies the following results:

\begin{Corollary}\label{betti splitting pd}
If $I = J + K$ is a Betti splitting ideal, then
\begin{enumerate}
\item $\reg(I) = \max\{\reg(J), \reg(K), \reg(J \cap K) - 1\}$,
\item $\mathrm{pdim}(I) = \max\{\mathrm{pdim}(J), \mathrm{pdim}(K), \mathrm{pdim}(J \cap K) + 1\}$.
\end{enumerate}
\end{Corollary}

\begin{Remark}\label{generally betti splitting}
Let $I, J,$ and $K$ be monomial ideals of $R$ such that $G(I)$ is the disjoint union of $G(J)$ and $G(K)$.
Furthermore, consider the following short exact sequence:
\begin{equation*}
(\ddagger)\hspace{1cm} 0 \rightarrow J \cap K \stackrel{\varphi}{\rightarrow} J \oplus K
\stackrel{\psi}{\rightarrow} J+K = I \rightarrow 0
\end{equation*}
where $\varphi(f) = (f,-f)$ and $\psi(g,h) = g+h$. The mapping cone construction applied to $(\ddagger)$ produces a free resolution of $I$. This implies the following inequality:
\[
\beta_i(I)  \leq \beta_i(J) + \beta_i(K) + \beta_{i-1}(J \cap K) \hspace{0.5cm} \text{for all } i \in \mathbb{N}.
\]
\end{Remark}

For a monomial $u \in R$, we denote by $\mathrm{supp}(u)$ the set of all variables dividing $u$. Furthermore, if a monomial ideal $I$ has minimal generating set $G(I) = \{u_1, \ldots, u_m\}$, we define $\mathrm{supp}(I) = \bigcup_{i=1}^m \mathrm{supp}(u_i)$.

The following lemma is well-known.
\begin{Lemma}\label{mergedlemma}
Let $I$ and $J$ be monomial ideals of $R$. We have
\begin{enumerate}
    \item If $\text{supp}(I) \cap \text{supp}(J) = \emptyset$, then $I \cap J = IJ$;
    \item For any monomial $u$ and integer $i$, $\beta_i(uI) = \beta_i(I)$.
\end{enumerate}
\end{Lemma}
\begin{Lemma}\label{splitting lemma}
Let $I$, $J$, and $K$ be monomial ideals of $R$ such that $\text{supp}(J) \cap \text{supp}(I) = \text{supp}(J) \cap \text{supp}(K) = \emptyset$, and $G(I) \cap G(JK) = \emptyset$. If every element $u \in G(I)$ can be written as $u = u_1u_2$ where $u_1 \in G(K)$, $\deg(u_2) > 0$, and $\text{supp}(u_2) \cap \text{supp}(K) = \emptyset$, then

\[
\beta_{i,\mathbf{a}}(I+JK) = \beta_{i,\mathbf{a}}(I)+\beta_{i,\mathbf{a}}(JK)+\beta_{i-1,\mathbf{a}}(JI) \hspace{.5cm} \text{for all } i\in \mathbb{N} \text{ and } \mathbf{a}\in\mathbb{N}^n;
\]

that is, $I + JK$ is a Betti splitting. In particularly $$\mathrm{pdim}(I+JK) = \max\{\mathrm{pdim}(I), \mathrm{pdim}(JK), \mathrm{pdim}(JI) + 1\}.$$
\end{Lemma}
\begin{proof}
Clearly, we have \(I\subseteq K\). Then, by Lemma \ref{mergedlemma}(1), we obtain
\[
I\cap JK = I\cap J\cap K = I\cap J = IJ.
\]
Then, by Corollary \ref{betti splitting pd}, we only need to prove that \(I + JK\) is a Betti splitting. Now, assume \(\beta_{i,\mathbf{a}}(JI)\neq 0\). Since \(\text{supp}(J) \cap \text{supp}(I)= \emptyset\), there exists \(x_i\in \text{supp}(J)\) such that \(x_i\) divides \(x^{\mathbf{a}}\), and thus \(\beta_{i,\mathbf{a}}(I)= 0\).

On the other hand, from \(\text{supp}(J) \cap \text{supp}(I)=\emptyset\), every minimal generator of \(IJ\) can be written as \(uv\), where \(u\in G(I)\) and \(v\in G(J)\). Further, it can be written as \(u_1u_2v\), where \(u_1\in G(K)\) and \(\text{supp}(u_2) \cap \text{supp}(K) = \emptyset\). Since \(\text{supp}(J) \cap \text{supp}(I)= \emptyset\), we have \(\text{supp}(u_2) \cap \text{supp}(J) = \emptyset\). Then there exists \(x_i\) that divides \(x^{\mathbf{a}}\) and does not belong to \(\text{supp}(J) \sqcup \text{supp}(K)\), so \(\beta_{i,\mathbf{a}}(JK)= 0\).

Finally, by using \cite[Theorem 2.3]{FHT}, we conclude that \(I + JK\) is a Betti splitting, as desired.
\end{proof}

\begin{Definition}\em\label{dominant de}
Let $I$ be a monomial ideal in $R$. We call $I$ a \emph{dominant ideal} if for any $m \in G(I)$ there exists a variable $x \in \text{supp}(m)$ such that for all $m' \in G(I) \backslash \{m\}$ we have $\deg_x m > \deg_x m'$, where for a monomial $v$ and a variable $x$, $\deg_x v$ denotes the exponent of $x$ in $v$.
\end{Definition}
\begin{Remark}\em
Note that a complete intersection monomial ideal is a dominant ideal, because its minimal generators do not share any common variables pairwise.
\end{Remark}

By \cite[Theorem 4.4]{A1}, if $I$ is a dominant ideal, then its Taylor resolution is minimal, leading to the following lemma.

\begin{Lemma}\em\label{dominant betti}
Let $I$ be a dominant ideal of $R$ and set $G(I) = \{u_1, \ldots, u_m\}$. Let the vector $\text{mdeg}(\text{lcm}\{u_1,\ldots,u_m\})$ be denoted by $\mathbf{b}$. Then $\text{pdim}(I) = m-1$ and its top multigraded Betti number is given by

\[
\beta_{m-1,\mathbf{a}}(I) = \left\{
\begin{array}{ll}
1, & \text{if } \mathbf{a} = \mathbf{b}; \\
0, & \text{otherwise.}
\end{array}
\right.
\]
In particular, when $i \neq m-1$,
\[
\beta_{i,\mathbf{b}}(I) = 0.
\]

Here, $\text{lcm}\{u_1,\ldots,u_m\}$ denotes the least common multiple of $u_1,\ldots,u_m$.

\end{Lemma}

\begin{Proposition} {\em (\cite[Proposition 3.1]{JKM})} \label{supp disjoint betti}\em
Let $S = \mathbb{K}[x_1, \dotsc,x_n, y_1\dotsc, y_m]$ be a polynomial ring. Let $I$ and $J$ be two monomial ideals, such that $\text{supp}(I)\cap \text{supp}(J)=\emptyset$. We might as well assume that $\text{supp}(I) \subseteq \{x_1, \ldots, x_n\}$ and $\text{supp}(J) \subseteq \{y_1, \ldots, y_m\}$.
\begin{enumerate}
\item[{\em(a)}] For $i \geq 1$, $\mathbf{a} \in \mathbb{N}^n$, $\mathbf{b} \in \mathbb{N}^m$ the following holds:
    \begin{align*}
        \beta_{i, (\mathbf{a},\mathbf{b})}^S (IJ) &= \sum_{\substack{j,k \geq 0 \\ j+k=i}} \beta_{j,\mathbf{a}}^S(I) \cdot \beta_{k,\mathbf{b}}^S(J) \\
        \mathrm{pdim}_S (IJ) &= \mathrm{pdim}_S (I) + \mathrm{pdim}_S (J)
    \end{align*}
\item[{\em(b)}] For $i \geq 0$ and $\mathbf{a} \in \mathbb{N}^n$, $\mathbf{b} \in \mathbb{N}^m$ the following holds:
    \begin{align*}
        \beta_{i, (\mathbf{a},\mathbf{b})}^S (I+J) &= \sum_{\substack{j,k \geq -1 \\ j+k=i-1}} \beta_{j,\mathbf{a}}^S(I) \cdot \beta_{k,\mathbf{b}}^S(J) \\
        \mathrm{pdim}_S (I+J) &= \mathrm{pdim}_S (I) + \mathrm{pdim}_S (J)+1.
    \end{align*}
\end{enumerate}
\end{Proposition}

The following lemma is important because it serves as the cornerstone of the induced-subgraph approach. We provide a proof for the completeness.

\begin{Lemma}\label{start}\em Let $H$ be  an induced subgraph of a weighted oriented graph $D$.
Then, for $\mathbf{a} \in \mathbb{N}^{V(D)}$ with $\mathrm{supp}(\mathbf{a}) \subseteq V(H)$, where $\mathrm{supp}(\mathbf{a})$ is defined as $\mathrm{supp}(x^{\mathbf{a}})$, we have
\[
\beta_{i,\mathbf{a}}(I(H)) = \beta_{i,\mathbf{a}}(I(D)) \quad \text{for all } i \geq 0.
\]
\end{Lemma}
\begin{proof}
Denote the monomial $\Theta = \prod_{x_i \in V(H)} x_i^{w(x_i)}$. Then, we have $I(H) = (I(D)_{\leq \Theta})$. Thus, the conclusion can be directly obtained from \cite[Proposition 56.1]{P}.
\end{proof}

\begin{Conventions}\label{conventions}
 We identify a monomial  with its multi-degree.  For example, if $\mathbf{a},\mathbf{b}\in\mathbb{N}^{V(H)}$ we will  write $\beta_{i,x^{\mathbf{a}}x^{\mathbf{b}}}(I(H))$ for $\beta_{i, \mathbf{a}+\mathbf{b}}(I(H))$.
\end{Conventions}

\begin{Definition}\em Let $H$ be a weighted oriented graph. We define the monomial
 $$
\Theta_H = \mathrm{lcm}\{u \mid u \in G(I(H))\}.
 $$
\end{Definition}

\section{Upper Bound}

In this section, we recursively provide an upper bound for the total Betti numbers of the edge ideal of a weighted oriented graph $G_n$. In the next section, we will prove that this upper bound is exact.

Recall that the edge ideal $I_n$ of $G_n$ is defined as
\[ I_n = (x_i y_j^{w_j} \mid 1 \leq i, j \leq n, i \neq j). \]

\begin{Lemma} \label{crowncolon}
We denote $I_n = I_{n-1} + (x_n)A + (y_n^{w_n})B$, where $A = (y_1^{w_1}, \ldots, y_{n-1}^{w_{n-1}})$ and $B = (x_1, \ldots, x_{n-1})$ are two monomial ideals.
For any $1 \leq s \leq n$, let $Q_s = (I_{n-1}, x_nx_1y_1^{w_1}, \ldots, x_nx_sy_s^{w_s})$ and $P_s = Q_{s-1} : x_nx_sy_s^{w_s}$. Note that $Q_0 = I_{n-1}$. Then the following statements hold:
\begin{enumerate}
  \item [$(1)$] $(I_{n-1} + x_nA) \cap y_n^{w_n}B = y_n^{w_n}C$, where $C = (I_{n-1}, x_nx_1y_1^{w_1}, \ldots, x_nx_{n-1}y_{n-1}^{w_{n-1}})$;
  \item [$(2)$] $P_s = (x_1, \ldots, x_{s-1}, x_{s+1}, \ldots, x_{n-1}, y_1^{w_1}, \ldots, y_{s-1}^{w_{s-1}}, y_{s+1}^{w_{s+1}}, \ldots, y_{n-1}^{w_{n-1}})$.
\end{enumerate}
\end{Lemma}
\begin{proof}
(1) Consider the following sequence of equalities:
\begin{align*}
(I_{n-1} + (x_n)A) \cap (y_n^{w_n})B
&= I_{n-1} \cap (y_n^{w_n})B + (x_n)A \cap (y_n^{w_n})B \\
&= (y_n^{w_n})I_{n-1} + (x_ny_n^{w_n})AB \quad \text{(since $I_{n-1} \subseteq B$ and Lemma~\ref{mergedlemma})}(1) \\
&= y_n^{w_n}\big(I_{n-1} + (x_n)AB\big)
\end{align*}

We observe that the generating set satisfies:
\[
G\big((x_n)AB\big) = \{x_nx_i y_j^{w_j} \mid 1 \leq i, j \leq n-1\}
\]
where the condition $x_nx_i y_j^{w_j} \in I_{n-1}$ holds whenever $i \neq j$. This leads to the final expression:
\[
(I_{n-1} + (x_n)A) \cap (y_n^{w_n})B = y_n^{w_n}\big(I_{n-1}, x_nx_1y_1^{w_1}, \dots, x_nx_{n-1}y_{n-1}^{w_{n-1}}\big)
\]
as desired.

(2) The assertion follows immediately from the fact that $Q_{s-1} : x_nx_sy_s^{w_s}$ is generated by the monomials $u_i/ \gcd(u_i,x_nx_sy_s^{w_s})$, for all $u_i\in G(Q_{s-1})$, see \cite[Proposition 1.2.2.]{HH}. Here, $\gcd(u_i, x_nx_sy_s^{w_s})$ denotes the greatest common divisor of $u_i$ and $x_nx_sy_s^{w_s}$ in the polynomial ring.
\end{proof}

\begin{Lemma} \label{binom}$m\geq0$ $n\geq2$ Then, we have
\[
\binom{2n-4}{m} = \sum_{k=2}^n 2^{m+4-2k}\binom{n-k}{m+4-2k}\binom{n-2}{k-2}.
\]
    \end{Lemma}

\begin{proof}This is equivalent to proving that for any $m \geq 0$ and $n \geq 0$,
\[
\binom{2n}{m} = \sum_{k=0}^{n} 2^{m-2k} \binom{n-k}{m-2k} \binom{n}{k}.
\]
We consider the number of ways to choose $m$ elements from $2n$ elements. Partition the $2n$ elements into $n$ disjoint pairs: $\{x_1, y_1\}, \ldots, \{x_n, y_n\}$. We categorize the selections based on the number of complete pairs $k$ included in the chosen subset:

1. \textbf{Select $k$ complete pairs}: This can be done in $\binom{n}{k}$ ways.

2. \textbf{Choose remaining $m-2k$ elements}: From the remaining $n-k$ pairs, select $m-2k$ pairs and independently pick one element from each selected pair ($x_j$ or $y_j$). This contributes $\binom{n-k}{m-2k} \cdot 2^{m-2k}$ possibilities.

Summing over all valid values of $k$ gives the total number of combinations:
\[
\binom{2n}{m} = \sum_{k=0}^{n} \binom{n}{k} \cdot 2^{m-2k} \cdot \binom{n - k}{m - 2k}.
\]
This completes the proof.\end{proof}
We  are  now ready to present the main result of this section.
\begin{Theorem} \label{main1}
For all $0\leq i$, we have
\[
\beta_i(I_n) \leq \sum_{k=2}^n (k-1)2^{i+3-2k}\binom{n}{k}\binom{n-k}{i+3-2k}+(2^{i+2}-2)\binom{n}{i+2}.
\]
\end{Theorem}
\begin{proof}
We proceed by induction on both $i$ and $n$. When $i = 0$, the right-hand side of the inequality is exactly the number of minimal generators of $I_n$, which is  $n(n-1)$. When $n = 2$, since the ideal $I_2$ is a complete intersection, it admits a minimal Taylor resolution. Consequently, it can be straightforwardly verified that the inequality holds.

     Now, assume $i>0$ and  $n>2$. We have
  \begin{align*}\tag{$\Delta$}
      \beta_i(I_n)
       &\leq\beta_i(I_{n-1}+x_nA)+\beta_i(y_n^{w_n}B)+\beta_{i-1}(y_n^{w_n}C)\\
      \qquad
      &\leq\beta_i(I_{n-1})+\beta_i(x_nA)+\beta_{i-1}(x_nI_{n-1})+\beta_i(y_n^{w_n}B)+\beta_{i-1}(y_n^{w_n}C)\\
      \qquad
      &=\beta_i(I_{n-1})+\beta_i(A)+\beta_{i-1}(I_{n-1})+\beta_i(B)+\beta_{i-1}(C).
    \end{align*}
The first inequality follows from Lemma~\ref{crowncolon}(1) and Remark~\ref{generally betti splitting}, the second inequality is obtained via $I_{n-1} \cap x_nA = x_nI_{n-1}$ combined with Remark~\ref{generally betti splitting}, and the final equality uses Lemma~\ref{mergedlemma}(2).

Now, consider the following short exact sequences:
$$
0\rightarrow R/P_s[-(w_s+2)]\rightarrow R/Q_{s-1}\rightarrow R/Q_s\rightarrow 0, \quad \text{for all}\ s=1,\ldots,n-1.
$$
Utilizing the mapping cone construction (see e.g. \cite[section 27]{P} for the details), we can deduce that
$$
\beta_{i}(R/Q_s)\leq \beta_{i-1}(R/P_s)+\beta_{i}(R/Q_{s-1}).
$$
Repeatedly applying the above inequality, we arrive at
$$
\beta_i(R/Q_{n-1})\leq \beta_{i-1}(R/P_{n-1})+\cdots+\beta_{i-1}(R/P_1)+\beta_{i}(R/Q_0)
$$
and
$$
\beta_{i-1}(C)\leq \beta_{i-2}(P_{n-1})+\cdots+\beta_{i-2}(P_1)+\beta_{i-1}(I_{n-1}).
$$
This implies that, for all non-negative integers $i$, the following inequality holds, which will be crucial in subsequent discussions:
\begin{align*}
  \beta_{i-1}(C)
  &\leq (n-1)\binom{2n-4}{i-1} + \beta_{i-1}(I_{n-1}) \\
  &\leq (n-1)\sum_{k=2}^n \left[2^{i+3-2k} \binom{n-2}{k-2} \binom{n-k}{i+3-2k}\right] + \beta_{i-1}(I_{n-1}) \\
  &\quad \text{(by Lemma \ref{binom})} \\
  &\leq \sum_{k=2}^n (k-1)2^{i+3-2k} \binom{n-1}{k-1} \binom{n-k}{i+3-2k} + \beta_{i-1}(I_{n-1}) \\
  &\quad \text{(using the combinatorial identity $p\binom{q}{p} = q\binom{q-1}{p-1}$)}
\end{align*}

Combining this with $(\Delta)$, we have
\allowdisplaybreaks
\begin{align*}
    \beta_i(I_n)
    &\leq \beta_i(I_{n-1}) + 2\beta_{i-1}(I_{n-1}) + 2\binom{n-1}{i+1} + (n-1)\binom{2n-4}{i-1} \\
    &\leq \sum_{k=2}^{n-1} (k-1)2^{i+3-2k}\binom{n-1}{k}\binom{n-1-k}{i+3-2k} + (2^{i+2}-2)\binom{n-1}{i+2} \\
    &\quad + 2\sum_{k=2}^{n-1} (k-1)2^{i+2-2k}\binom{n-1}{k}\binom{n-1-k}{i+2-2k} + 2(2^{i+1}-2)\binom{n-1}{i+1} \\
    &\quad + 2\binom{n-1}{i+1} + \sum_{k=2}^n (k-1)2^{i+3-2k}\binom{n-1}{k-1}\binom{n-k}{i+3-2k} \\
    &= \sum_{k=2}^{n-1} (k-1)2^{i+3-2k}\binom{n-1}{k}\binom{n-k}{i+3-2k} + (2^{i+2}-2)\binom{n}{i+2} \\
    &\quad + \sum_{k=2}^n (k-1)2^{i+3-2k}\binom{n-1}{k-1}\binom{n-k}{i+3-2k} \\
    &= \sum_{k=2}^{n-1} (k-1)2^{i+3-2k}\binom{n}{k}\binom{n-k}{i+3-2k} + (n-1)2^{i+3-2n}\binom{n-1}{n-1}\binom{n-n}{i+3-2n} \\
    &\quad + (2^{i+2}-2)\binom{n}{i+2} \\
    &= \sum_{k=2}^{n} (k-1)2^{i+3-2k}\binom{n}{k}\binom{n-k}{i+3-2k} + (2^{i+2}-2)\binom{n}{i+2}.
\end{align*}
The second inequality above uses the induction hypothesis, and the first two equalities apply the combinatorial identity $\binom{q-1}{p} + \binom{q-1}{p-1} = \binom{q}{p}$. This completes the proof by induction.
\end{proof}

\section{Betti Numbers of Edge Ideals of Induced Subgraphs}

In this section, we are dedicated to computing the projective dimension and the top multigraded Betti numbers of the edge ideals of some induced subgraphs that may be derived from the weighted crown graph $G_n$. These graphs include the weighted crown graph $G_s$, graphs of type $G_{s,t}$, graphs of type $G_{s,t}'$, graphs of type $G_{m,s,t}$, and so on.

\subsection{Induced Subgraphs}

\begin{itemize}
    \item[Type 1:] Weighted oriented crown graph $G_s$

    Recall that the edge ideal $I_n$ of $G_n$ is defined as
\[ I_n = (x_i y_j^{w_j} \mid 1 \leq i, j \leq n, i \neq j). \]

Below, (1) of the following lemma can be verified directly using \cite[Proposition 1.2.2.]{HH}, and (2) can be obtained directly from the definition of $I_s$.
\begin{Lemma}\label{lemma}
 For all $s\geq 3$, we have
 \begin{enumerate}
   \item $I_s:x_s=(y_1^{w_1},\ldots,y_{s-1}^{w_{s-1}}, x_1y_s^{w_s}, \ldots, x_{s-1}y_s^{w_s} )=A$;
   \vspace{2mm}

   \item $(I_s, x_s)=(I_{s-1}, y_s^{w_s}(x_1,\ldots,x_{s-1}), x_s)=B$.
 \end{enumerate}
              \end{Lemma}
  Obviously, $\lcm(y_1^{w_1}, \ldots, y_{s-1}^{w_{s-1}}, x_1y_s^{w_s}, \ldots, x_{s-1}y_s^{w_s})=
y_1^{w_1} \ldots y_{s-1}^{w_{s-1}} y_s^{w_s} x_1 \ldots x_{s-1}$,
which we denote as $u$. Note that $\Theta_{G_s} = ux_s$.

\begin{Proposition}\label{Gs betti}
Then $\mathrm{pdim}(I_s)=2s-3$ and its top multigraded Betti number is given by
\[
\beta_{2s-3,\alpha}(I_s) = \left\{
\begin{array}{ll}
s-1, & \text{if } \alpha = \Theta_{G_s} ; \\
0, & \text{otherwise.}
\end{array}
\right.
\]
\end{Proposition}
   \begin{proof} We proceed by induction on $s$. When $s=2$, since $I_2 = (x_1y_2^{w_2}, x_2y_1^{w_1})$ is a complete intersection, the conclusion follows from Lemma \ref{dominant betti}. Next, assume $s > 2$. By the inductive hypothesis, we have $\text{pdim}(I_{s-1}) = 2s-5$ and
\[
\beta_{2s-5,\alpha}(I_{s-1}) = \left\{
\begin{array}{ll}
s-2, & \text{if } \alpha = \Theta_{G_{s-1}}; \\
0, & \text{otherwise.}
\end{array}
\right.
\]

Since $\mathrm{supp}(y_s^{w_s}) \cap \mathrm{supp}(x_1,\ldots,x_{s-1}, x_s) = \mathrm{supp}(y_s^{w_s}) \cap \mathrm{supp}(I_{s-1}) = \emptyset$ and both $(y_s^{w_s})$ and $(x_1,\ldots,x_{s-1}, x_s)$ are dominant ideals, we can apply Proposition \ref{supp disjoint betti} and Lemma \ref{dominant betti} to obtain $\mathrm{pdim}(y_s^{w_s}(x_1,\ldots,x_{s-1})) = s-2$ and $\mathrm{pdim}(y_s^{w_s}I_{s-1}) = 2s-5$. Furthermore,
\[
\beta_{2s-5,\alpha}(y_s^{w_s}I_{s-1}) = \left\{
\begin{array}{ll}
s-2, & \text{if } \alpha = u; \\
0, & \text{otherwise.}
\end{array}
\right.
\]

Let $I_{s-1}$, $(y_s^{w_s})$, and $(x_1,\ldots,x_{s-1})$ be $I$, $J$, and $K$ respectively in Lemma \ref{splitting lemma}. It can be verified that $I + JK$ satisfies the conditions of Lemma \ref{splitting lemma}. Therefore, we have
\begin{align*}
\mathrm{pdim}\big(I_{s-1} + y_s^{w_s}(x_1,\ldots,x_{s-1})\big)
&= \max\{2s-5,\ s-2,\ 2s-5+1\} \\
&= 2s-4
\end{align*}
and
\begin{align*}
\beta_{2s-4,\alpha}\big(I_{s-1} + y_s^{w_s}(x_1,\ldots,x_{s-1})\big)
&= \beta_{2s-5,\alpha}(y_s^{w_s}I_{s-1}) \\
&= \begin{cases}
    s-2, & \text{if } \alpha = u; \\
    0,   & \text{otherwise.}
   \end{cases}
\end{align*}

By applying part (b) of Proposition \ref{supp disjoint betti} to $(I_{s-1} + y_s^{w_s}(x_1,\ldots,x_{s-1})) + (x_s)$, we obtain
\[
\mathrm{pdim}(R/B) = 2s-2
\]
and
\[
\beta_{2s-2,\alpha}(R/B) = \left\{
\begin{array}{ll}
s-2, & \text{if } \alpha = \Theta_{G_s}; \\
0, & \text{otherwise.}
\end{array}
\right.
\]

According to Definition \ref{dominant de}, $A$ is a dominant ideal. Therefore, by applying Lemma \ref{dominant betti}, we obtain
\[
\mathrm{pdim}(R/A) = 2s-2
\]
and
\[
\beta_{2s-2,\alpha}(R/A) = \left\{
\begin{array}{ll}
1, & \text{if } \alpha = u; \\
0, & \text{otherwise.}
\end{array}
\right.
\]
According to \cite[lemma~2.1]{HT}, we have
\[
\mathrm{pdim}(R/I_s) = 2s-2.
\]
Furthermore, we consider the following short exact sequence of multigraded $R$-modules:
\begin{equation*}\label{first}
0 \longrightarrow \frac{R}{I_s:x_s}[-x_s] \longrightarrow \frac{R}{I_s} \longrightarrow \frac{R}{(I_s,x_s)} \longrightarrow 0.
\end{equation*}
By Lemma \ref{lemma}, this is equivalent to
\begin{equation*}\label{second}
0 \longrightarrow \frac{R}{A}[-x_s] \longrightarrow \frac{R}{I_s} \longrightarrow \frac{R}{B} \longrightarrow 0.
\end{equation*}

The above short exact sequence induces the long exact sequence:
\begin{equation*}
\begin{aligned}
&0 \longrightarrow \Tor_{2s-2}^R(R/A[-x_s],\mathbb{K})_{\alpha} \longrightarrow \Tor_{2s-2}^R(R/I_s,\mathbb{K})_{\alpha} \\
&\longrightarrow \Tor_{2s-2}^R(R/B,\mathbb{K})_{\alpha} \longrightarrow \Tor_{2s-3}^R(R/A[-x_s],\mathbb{K})_{\alpha} \longrightarrow \ldots
\end{aligned}
\end{equation*}

Obviously, for all $\alpha \neq \Theta_{G_s}$, we have $\Tor_{2s-2}^R(R/I_s,\mathbb{K})_{\alpha} = 0$.

Since $A$ is a dominant ideal, by Lemma~\ref{dominant betti}, we have \[\Tor_{2s-3}^R(R/A[-x_s],\mathbb{K})_{\Theta_{G_s}} = 0.\]
Thus, we obtain
\begin{align*}
\beta_{2s-2,\Theta_{G_s}}(R/I_s)
&= \beta_{2s-2,\Theta_{G_s}}(R/A[-x_s]) + \beta_{2s-2,\Theta_{G_s}}(R/B) \\
&= (s-2) + 1 \\
&= s-1
\end{align*}
which completes the proof. \end{proof}
    \item[Type 2:] Weighted oriented unbalanced crown graph $G_{s,t}$

    For any integers $1 < s \neq t$, an \emph{$(s,t)$-weighted oriented unbalanced crown graph} is defined as $G_{s,t} = (V, E, w)$, where:
\begin{itemize}
    \item $V(G_{s,t}) = \{x_1,\ldots,x_s\} \cup \{y_1,\ldots,y_t\}$ is the vertex set;
    \item $E(G_{s,t}) = \{(x_i,y_j) \mid 1 \leq i \leq s,\ 1 \leq j \leq t,\ i \neq j\}$ is the edge set;
    \item $w: V(G_{s,t}) \to \mathbb{N}^+$ is the weight function.
\end{itemize} In \cite{R}, the underlying graph was referred to as an \emph{unbalanced crown graph}, and its graded Betti numbers for the edge ideal were explicitly determined.

We assume $1 < t < s$.
The edge ideal of the graph $G_{s,t}$ is denoted by $I_{s,t}$ and is given by
\[
I_{s,t} = (x_i y_j^{w_j} \mid 1 \leq i \leq s, 1 \leq j \leq t, i \neq j).
\]
Therefore, we have $\Theta_{G_{s,t}}=x_1\ldots x_sy_1^{w_1}\ldots y_t^{w_t}$.
\begin{Proposition}\label{Gst betti}
Then $\mathrm{pdim}(I_{s,t})=s+t-3$ and its top multigraded Betti number is given by
\[
\beta_{s+t-3,\alpha}(I_{s,t}) = \left\{
\begin{array}{ll}
t-1, & \text{if } \alpha = \Theta_{G_{s,t}} ; \\
0, & \text{otherwise.}
\end{array}
\right.
\]
\end{Proposition}
\begin{proof}Denote $B = (y_1^{w_1}, \ldots, y_t^{w_t})$  and  $C = (x_{t+1}, \ldots, x_s)$. Then, from the definition of $I_{s,t}$,  we have $I_{s,t} = I_t + CB$.
According to Proposition \ref{Gs betti}, the projective dimension of $I_t$ is $2t-3$, and its top multigraded Betti number is:
\[
\beta_{2t-3,\alpha}(I_t) =
\begin{cases}
t-1, & \text{if } \alpha = \Theta_{G_t}, \\
0, & \text{otherwise.}
\end{cases}
\]
Since $\mathrm{supp}(C) \cap \mathrm{supp}(I_t) = \mathrm{supp}(C) \cap \mathrm{supp}(B) = \emptyset$ and both $C$ and $B$ are dominant ideals, we can apply Proposition \ref{supp disjoint betti} and Lemma \ref{dominant betti} to obtain $\mathrm{pdim}(CB) = s-2$ and $\mathrm{pdim}(CI_t) = s+t-4$. Furthermore,
\[
\beta_{s+t-4,\alpha}(CI_t) = \left\{
\begin{array}{ll}
t-1, & \text{if } \alpha = \Theta_{G_{s,t}}; \\
0, & \text{otherwise.}
\end{array}
\right.
\]

Let $I_t$, $C$, and $B$ be $I$, $J$, and $K$ respectively in Lemma \ref{splitting lemma}. It can be verified that $I + JK$ satisfies the conditions of Lemma \ref{splitting lemma}. Therefore, we have
\[
\mathrm{pdim}(I_{s,t}) = \max\{2t-3, s-2, s+t-4+1\}=s+t-3
\]
and
\[
\beta_{s+t-3,\alpha}(I_{s,t}) = \beta_{s+t-4,\alpha}(CI_t)= \left\{
\begin{array}{ll}
t-1, & \text{if } \alpha = \Theta_{G_{s,t}}; \\
0, & \text{otherwise.}
\end{array}
\right.
\]

\end{proof}
    \item[Type 3:] Weighted oriented generalized crown graph $G_{m,s,t}$

    For any integers $1 < m <s, t$, an $(m,s,t)$-\emph{Weighted oriented generalized crown graph} is defined as $G_{m,s,t} = (V(G_{m,s,t}), E(G_{m,s,t}), w)$, where:
\begin{itemize}
    \item $V(G_{m,s,t}) = \{x_1, \ldots, x_s, y_1, \ldots, y_t\}$ is the vertex set;
    \item $E(G_{m,s,t}) = \Big\{(x_i, y_j) \mid 1 \leq i \leq s, 1 \leq j \leq t, i \neq j \Big\}\cup \Big\{(x_i, y_j) \mid m+1 \leq i \leq s, 1 \leq j \leq t\Big\}$ is the edge set;
    \item $w: V(G_{m,s,t}) \to \mathbb{N}^+$ is the weight function.
\end{itemize}
The edge ideal of the graph $G_{m,s,t}$ is denoted by $I_{m,s,t}$ and is given by
\[
\begin{split}
I_{m,s,t} &= (x_i y_j^{w_j} \mid 1 \leq i \leq s, 1 \leq j \leq t, i \neq j) \\
&\quad + (x_i y_j^{w_j} \mid m+1 \leq i \leq s, 1 \leq j \leq t).
\end{split}
\]
Therefore, we have $\Theta_{G_{m,s,t}}=x_1\ldots x_sy_1^{w_1}\ldots y_t^{w_t}$.
\begin{Proposition}\label{Gmst betti}
Then $\mathrm{pdim}(I_{m,s,t})=s+t-3$ and its top multigraded Betti number is given by
\[
\beta_{s+t-3,\alpha}(I_{m,s,t}) = \left\{
\begin{array}{ll}
m-1, & \text{if } \alpha = \Theta_{G_{m,s,t}} ; \\
0, & \text{otherwise.}
\end{array}
\right.
\]
\end{Proposition}

\begin{proof}Denote $A = (x_1, \ldots, x_s)$, and  $B= (y_{m+1}^{w_{m+1}}, \ldots, y_t^{w_t})$. Then, from the definition of $I_{m,s,t}$,  we have $I_{m,s,t} =I_{s,m} + BA$.
According to Proposition \ref{Gs betti}, the projective dimension of $I_{s,m}$ is $s+m-3$, and its top multigraded Betti number is:
\[
\beta_{s+m-3,\alpha}(I_{s,m}) =
\begin{cases}
m-1, & \text{if } \alpha = \Theta_{G_{s,m}}, \\
0, & \text{otherwise.}
\end{cases}
\]
Since $\mathrm{supp}(B) \cap \mathrm{supp}(I_{s,m}) = \mathrm{supp}(B) \cap \mathrm{supp}(A) = \emptyset$ and both $B$ and $A$ are dominant ideals, we can apply Proposition \ref{supp disjoint betti} and Lemma \ref{dominant betti} to obtain $\mathrm{pdim}(BA) = s+t-m-2$ and $\mathrm{pdim}(BI_{s,m}) = s+t-4$. Furthermore,
\[
\beta_{s+t-4,\alpha}(BI_{s,m}) = \left\{
\begin{array}{ll}
m-1, & \text{if } \alpha = \Theta_{G_{m,s,t}}; \\
0, & \text{otherwise.}
\end{array}
\right.
\]

Let $I_{s,m}$, $B$, and $A$ be $I$, $J$, and $K$ respectively in Lemma \ref{splitting lemma}. It can be verified that $I + JK$ satisfies the conditions of Lemma \ref{splitting lemma}. Therefore, we have
\begin{align*}
\mathrm{pdim}(I_{m,s,t})
&= \max\{s+m-3,\ s+t-m-2,\ s+t-4+1\} \\
&= s+t-3
\end{align*}
and
\[
\beta_{s+t-3,\alpha}(I_{m,s,t}) = \beta_{s+t-4,\alpha}(BI_{s,m})= \left\{
\begin{array}{ll}
t-1, & \text{if } \alpha = \Theta_{G_{m,s,t}}; \\
0, & \text{otherwise.}
\end{array}
\right.
\]\end{proof}

    \item[Type 4:] Weighted oriented complete bipartite graph $K_{s,t}$

      For any integers $1 \leq s, t$, an $(s,t)$-\emph{Weighted oriented complete bipartite graph} is defined as $K_{s,t} = (V(K_{s,t}), E(K_{s,t}), w)$, where:
\begin{itemize}
    \item $V(K_{s,t}) = \{x_1, \ldots, x_s, y_1, \ldots, y_t\}$ is the vertex set;
    \item $E(K_{s,t}) = \Big\{(x_i, y_j) \mid 1 \leq i \leq s, 1 \leq j \leq t\Big\}$ is the edge set;
    \item $w: V(K_{s,t}) \to \mathbb{N}^+$ is the weight function.
\end{itemize}

The edge ideal of the graph $K_{s,t}$ is denoted by $J_{s,t}$ and is given by
\[
J_{s,t} = (x_i y_j^{w_j} \mid 1 \leq i \leq s, 1 \leq j \leq t).
\]
Therefore, we have $\Theta_{K_{s,t}}=x_1\ldots x_sy_1^{w_1}\ldots y_t^{w_t}$.

\begin{Proposition}\label{Kst betti}
Then $\mathrm{pdim}(J_{s,t})=s+t-2$ and its top multigraded Betti number is given by
\[
\beta_{s+t-2,\alpha}(J_{s,t}) = \left\{
\begin{array}{ll}
1, & \text{if } \alpha = \Theta_{K_{s,t}} ; \\
0, & \text{otherwise.}
\end{array}
\right.
\]
\end{Proposition}
\begin{proof}
Denote $A = (x_1, \ldots, x_s)$, and  $B= (y_1^{w_1}, \ldots, y_t^{w_t})$. Then, from the definition of $J_{s,t}$,  we have $J_{s,t} =AB$.
Since $A$ and $B$ are dominant ideals, we can apply Proposition \ref{supp disjoint betti} and Lemma \ref{dominant betti} to obtain  $\mathrm{pdim}(J_{s,t})=\mathrm{pdim}(AB) = s+t-2$. Furthermore,
\[
\beta_{s+t-2,\alpha}(J_{s,t})=\beta_{s+t-2,\alpha}(AB) = \left\{
\begin{array}{ll}
1, & \text{if } \alpha = \Theta_{K_{s,t}} ; \\
0, & \text{otherwise.}
\end{array}
\right.
\]
\end{proof}
\end{itemize}

\begin{Notation}
We define:
\begin{itemize}
    \item $\mathcal{H}_1$ - the collection of induced subgraphs of $G_n$ of the first three types;
    \item $\mathcal{H}_2$ - the collection of induced subgraphs of $G_n$ of the fourth type.
\end{itemize}
\end{Notation}

Combining Propositions~\ref{Gs betti}, \ref{Gst betti}, \ref{Gmst betti}, and~\ref{Kst betti} yields:
\begin{Proposition}\label{inducetop}\em
Let $H$ be an induced subgraph of $G_n$.
\begin{enumerate}
    \item If $H \in \mathcal{H}_1$ contains exactly $k$ complete pairs of the form $\{x_i,y_i\}$ in $V(H)$, then
    $\mathrm{pdim}(I(H)) = |V(H)| - 3$ and its top multigraded Betti number is given by
    \[
    \beta_{|V(H)|-3, \alpha}(I(H)) = \begin{cases}
        k-1, & \text{if } \alpha = \Theta_{H}; \\
        0, & \text{otherwise.}
    \end{cases}
    \]

    \item If $H \in \mathcal{H}_2$, then $\mathrm{pdim}(I(H)) = |V(H)| - 2$ and its top multigraded Betti number is given by
    \[
    \beta_{|V(H)|-2, \alpha}(I(H)) = \begin{cases}
        1, & \text{if } \alpha = \Theta_{H}; \\
        0, & \text{otherwise.}
    \end{cases}
    \]
\end{enumerate}
\end{Proposition}

\subsection{Main results}

\begin{Definition}\label{topsupp}
For each integer $2 \leq k \leq n$, define:
\[
\mathcal{N}_i^k(G_n) := \left\{\Theta_H \middle|
\begin{aligned}
& H \in \mathcal{H}_1,\ V(H)\text{ contains exactly }k \\
& \text{complete pairs }\{x_i,y_i\},\ |V(H)|=i+3
\end{aligned}
\right\},
\]
and let $\mathcal{N}_i(G_n) := \bigcup_{k=2}^n \mathcal{N}_i^k(G_n)$.

For $\mathcal{H}_2$-type subgraphs, define:
\[
\mathcal{M}_i(G_n) := \{\Theta_H \mid H \in \mathcal{H}_2,\ |V(H)|=i+2\}.
\]
\end{Definition}

\begin{Theorem}\label{main2}
For all $0 \leq i \leq 2n-3$, the $i$-th total Betti number is:
\[
\beta_i(I_n) = \sum_{k=2}^n (k-1)2^{i+3-2k}\binom{n}{k}\binom{n-k}{i+3-2k} + (2^{i+2}-2)\binom{n}{i+2}.
\]
Moreover, the total Betti numbers of the weighted crown graph $G_n$ are independent of the vertex weights.
\end{Theorem}

\begin{proof}
The base case $i=0$ was established in Theorem~\ref{main1}. Fix $1 \leq i \leq 2n-3$.

Following the proof technique of Lemma~\ref{binom}, we partition the vertex set as
$V(G_n) = \bigsqcup_{j=1}^n \{x_j,y_j\}$. For any subset $T \subseteq V(G_n)$, our classification yields:

\begin{enumerate}
\item From our classification of the first three graph types, we have the following characterization:
\[
G_n[T] \in \mathcal{H}_1 \iff \exists k \geq 2 \text{ such that } T \text{ contains exactly } k \text{ pairs } \{x_i, y_i\}.
\]

By Proposition~\ref{inducetop}, any subgraph $H = G_n[T] \in \mathcal{H}_1$ with $\mathrm{pdim}(I(H)) = i$ satisfies $|V(H)| = i + 3$.
Applying Lemma~\ref{binom}'s counting argument, for each $2 \leq k \leq n$ there exist precisely
\[
a_i^k := 2^{i+3-2k}\binom{n}{k}\binom{n-k}{i+3-2k}
\]
induced subgraphs $\{H_j^k\}_{j=1}^{a_i^k}$ in $\mathcal{H}_1$ having $\mathrm{pdim}(I(H_j^k)) = i$. Consequently, we have the cardinality $|\mathcal{N}_i^k(G_n)| = a_i^k$.   By Proposition~\ref{inducetop}(1) and Lemma~\ref{start}, we obtain the Betti number equality:
\[
\beta_{i,\Theta_{H_{\ell}^k}}(I_n) = \beta_{i,\Theta_{H_{\ell}^k}}(I(H_{\ell}^k)) = k-1
\]
for any $\Theta_{H_{\ell}^k} \in \mathcal{N}_i^k(G_n)$ and $2 \leq k \leq n$.

This immediately yields the following summation:
\begin{align*}\tag{$\dag$}
\sum_{\alpha\in \mathcal{N}_i(G_n)}\beta_{i,\alpha}(I_n)
&= \sum_{k=2}^n\sum_{\alpha\in \mathcal{N}_i^k(G_n)}\beta_{i,\alpha}(I_n) \\
&= \sum_{k=2}^n\sum_{\alpha\in \mathcal{N}_i^k(G_n)}(k-1) \\
&= \sum_{k=2}^n a_i^k(k-1)
\end{align*}
where the final equality follows from the cardinality $|\mathcal{N}_i^k(G_n)| = a_i^k$.

\item From our classification of the fourth graph type, we establish the characterization:
\[
G_n[T] \in \mathcal{H}_2 \iff T \text{ contains no complete pair } \{x_i,y_i\}.
\]

By Proposition~\ref{inducetop}, any subgraph $H = G_n[T] \in \mathcal{H}_2$ with $\mathrm{pdim}(I(H)) = i$ satisfies $|V(H)| = i + 2$.
Applying Lemma~\ref{binom}'s counting argument, there exist precisely
\[
b_i := (2^{i+2} - 2)\binom{n}{i+2}
\]
induced subgraphs $\{H'_j\}_{j=1}^{b_i}$ in $\mathcal{H}_2$ having $\mathrm{pdim}(I(H'_j)) = i$. Consequently, we have the cardinality $|\mathcal{M}_i(G_n)| = b_i$.

By Proposition~\ref{inducetop}(2) and Lemma~\ref{start}, we obtain the Betti number equality:
\[
\beta_{i,\Theta_{H'_{\ell}}}(I_n) = \beta_{i,\Theta_{H'_{\ell}}}(I(H'_{\ell})) = 1
\]
for any $\Theta_{H'_{\ell}} \in \mathcal{M}_i(G_n)$.

This immediately yields the following summation:
\begin{align*}\tag{$\ddag$}
\sum_{\alpha\in \mathcal{M}_i(G_n)}\beta_{i,\alpha}(I_n)
&= \sum_{\alpha\in \mathcal{M}_i(G_n)} 1 \\
&= |\mathcal{M}_i(G_n)| \\
&= b_i.
\end{align*}
\end{enumerate}
Observe that $\mathcal{N}_i(G_n) \cap \mathcal{M}_i(G_n) = \emptyset$. Combining equations ($\dag$) and ($\ddag$), we derive:
\begin{align*}
\beta_i(I_n) &\geq \sum_{\substack{\alpha \in \mathcal{N}_i(G_n)  \sqcup \mathcal{M}_i(G_n)}} \beta_{i,\alpha}(I_n) \\
&= \sum_{k=2}^n (k-1)a_i^k + b_i \\
&= \sum_{k=2}^n (k-1)2^{i+3-2k}\binom{n}{k}\binom{n-k}{i+3-2k} + (2^{i+2}-2)\binom{n}{i+2}.
\end{align*}
The theorem follows immediately from Theorem~\ref{main1}.
\end{proof}

From Theorem~\ref{main2}, we immediately obtain the main result of this paper:

\begin{Theorem}\label{main}
Let $G_n$ be the weighted oriented crown graph as described previously, and let $I_n$ be its edge ideal.  Then its multigraded Betti numbers are given by
\[
\beta_{i, \alpha}(I_n) = \begin{cases}
k-1, & \text{if } \alpha \in \mathcal{N}_i^k(G_n) \text{ for some }2\leq k\leq n; \\
1, & \text{if } \alpha \in \mathcal{M}_i(G_n); \\
0, & \text{otherwise.}
\end{cases}
\]
Furthermore, the graded Betti numbers of $I_n$ are
\[
\beta_{i,j}(I_n) =\sum_{k=2}^n \sum_{\substack{\alpha \in \mathcal{N}_i^k(G_n) \\ \deg(\alpha)=j}} (k-1) + |\{\alpha \in \mathcal{M}_i(G_n) : \deg(\alpha)=j\}|.
\]
\end{Theorem}

\begin{Corollary}\label{crown reg}
Let $G_n$ be the weighted oriented crown graph as described previously, and let $I_n$ be its edge ideal.  Then its regularity is
\[
\mathrm{reg}(I_n) = \sum_{i=1}^n w_i - n + 3.
\]
\end{Corollary}
\begin{proof}

From Proposition~\ref{Gs betti}, we have $\mathrm{pdim}(I_n) = 2n - 3$ and $\beta_{2n-3, \Theta_{G_n}}(I_n) = n - 1$.
Since
\[
\deg(\Theta_{G_n}) - (2n - 3) = \sum_{i=1}^n w_i + n - (2n - 3) = \sum_{i=1}^n w_i - n + 3,
\]
by Theorem~\ref{main}, it suffices to prove that for all $0 \leq i < 2n - 3$ and $\Theta_H \in \mathcal{N}_i(G_n) \cup \mathcal{M}_i(G_n)$,
\[
\deg(\Theta_H) - i \leq \deg(\Theta_{G_n}) - (2n - 3)
\]
holds. This inequality is equivalent to
\[
(2n - 3) - i \leq \deg(\Theta_{G_n}) - \deg(\Theta_H).
\]
Since $H$ is a proper induced subgraph of $G_n$ with $i + 2$ or $i + 3$ vertices, the above inequality must hold.
\end{proof}

From Theorem~\ref{main}, we recover the results of Rather and Singh in \cite{RS1,RS2} concerning the graded Betti numbers of edge ideals of crown graphs.

\begin{Corollary}{\em (\cite[Theorem 3.1]{RS1} and \cite[Theorem 3.5 and Lemma 3.6]{RS2})}\label{cor}
If every vertex of the weighted oriented crown graph $G_n$ has weight $1$, then the graded Betti numbers of $I_n$ are given by
\[
\beta_{i, j}(I_n) = \begin{cases}
(2^{i+2} - 2)\binom{n}{i+2}, & \text{if } j = i + 2; \\
\sum\limits_{k=2}^n (k-1)2^{i+3-2k}\binom{n}{k}\binom{n-k}{i+3-2k}, & \text{if } j = i + 3; \\
0, & \text{otherwise.}
\end{cases}
\]
\end{Corollary}
\begin{proof}
When all vertices have weight $1$, Definition~\ref{topsupp} implies that all monomials in $\mathcal{N}_i^k(G_n)$ (resp. $\mathcal{M}_i(G_n)$) have degree $i+3$ (resp. $i+2$). From Theorem~\ref{main}, we obtain:
\[
\beta_{i,j}(I_n) = \begin{cases}
|\mathcal{M}_i(G_n)|, & \text{if } j=i+2; \\
\sum_{k=2}^n (k-1)|\mathcal{N}_i^k(G_n)|, & \text{if } j=i+3; \\
0, & \text{otherwise.}
\end{cases}
\]
The proof of Theorem~\ref{main2} shows that $|\mathcal{M}_i(G_n)| = (2^{i+2}-2)\binom{n}{i+2}$ and $|\mathcal{N}_i^k(G_n)| = 2^{i+3-2k}\binom{n}{k}\binom{n-k}{i+3-2k}$, which completes the proof.
\end{proof}

\begin{Remark}
Let $H$ be a weighted oriented unbalanced crown graph, generalized crown graph, or complete bipartite graph as defined in Subsection~4.1. Since $H$ can be identified with an induced subgraph of $G_n$ via vertex relabeling, the multigraded Betti numbers of $I(H)$ are determined by Lemma~\ref{start} and Theorem~\ref{main}.
\end{Remark}

Our computations reveal that the multigraded Betti numbers of $I(G_n)$ are entirely determined by the top multigraded Betti numbers of its induced subgraphs, leading to:

\begin{Proposition}\label{po}
Let $G$ denote the weighted oriented crown graph $G_n$ studied in this paper. For any $i \geq 0$ and multidegree $\alpha$, the following implication holds:
\[
\tag{$\S$}
\beta_{i,\alpha}(I(G)) \neq 0 \implies \mathrm{pdim}\big(I(G[\mathrm{supp}(\alpha)])\big) = i.
\]
\end{Proposition}

This leads us to consider:
\begin{Question}
Can properties of induced subgraphs fully determine the minimal multigraded free resolution of $I(G_n)$?
\end{Question}

Motivated by the computational simplicity of edge ideals satisfying the implication ($\S$), we ask:
\begin{Question}
Characterization of weighted oriented graphs satisfying implication ($\S$) via combinatorial properties.
\end{Question}

{\bf \noindent Acknowledgment:}
This research is supported by NSFC (No. 11971338). We gratefully acknowledge the use of the computer algebra system Cocoa (\cite{CoCoA}) for our experiments. The first author would like to acknowledge Prof. Lijun Ji and Prof. Xinrong MA  for their guidance in the proof of combinatorial identities.

\vspace{2mm}

{\bf\noindent Statement:} On behalf of all authors, the corresponding author states that there is
no conflict of interest.

\end{document}